# Politiques de Tests Partiels & Systèmes de Sécurité

Évaluation et optimisation des probabilités de défaillance à la demande (PFD) des systèmes instrumentés de sécurité (SIS) soumis à des tests de révision partiels

**Florent Brissaud [1,2], Anne Barros [2] & Christophe Bérenguer [2]**

[1] *INERIS*
*Parc Technologique Alata, BP 2*
*60550 Verneuil-en-Halatte, France*
*florent.brissaud@ineris.fr*

[2] *UTT – ICD, FRE CNRS 2848*
*12 rue Marie Curie*
*10010 Troyes cedex, France*
*{prenom.nom}@utt.fr*

RÉSUMÉ. *Des expressions générales, obtenues par une approche analytique, sont proposées afin d'évaluer les probabilités de défaillance à la demande (PFD) des systèmes d'architecture MooN (i.e. k-out-of-n) soumis à des tests de révision partiels et complets. Les tests partiels (e.g. inspections visuels, essais imparfaits) ne permettent de détecter que certaines défaillances, tandis qu'à la suite des tests complets le système est rétabli dans des conditions « comme neuves ». En suivant la démarche proposée, une application présente un moyen d'estimer les performances du système et des politiques de tests à partir du retour d'expérience obtenu au cours des tests partiels et complets. Une optimisation de la répartition des tests partiels est également proposée, permettant de réduire la probabilité moyenne de défaillance du système à la demande (PFDavg).*

ABSTRACT. *A set of general formulas is proposed for the probability of failure on demand (PFD) assessment of MooN architecture (i.e. k-out-of-n) systems subject to partial and full tests. Partial tests (e.g. visual inspections, imperfect testing) may detect only some failures, whereas owing to a full test, the system is restored to an as good as new condition. Following the proposed approach and according to an example, performance estimations of the system and test policies are presented, by using the feedback from partial and full tests. An optimization of the partial test distribution is also proposed, which allows reducing the average probability of system failure on demand (PFDavg).*

MOTS-CLÉS : *Probabilité de défaillance à la demande, PFD, Système instrumenté de sécurité, SIS, Période de révision, Test partiel, Test complet.*

KEYWORDS: *Probability of failure on demand, PFD, Safety Instrumented System, SIS, Proof test interval, Partial test, Full test.*

**1. Introduction**

Les systèmes instrumentés de sécurité (SIS) jouent un rôle clef dans la prévention des risques industriels, en tant que barrières de sécurité. L'objectif d'un SIS est d'assurer une ou plusieurs fonctions destinées à assurer ou à maintenir un état sûr d'un équipement commandé (EUC), par rapport à un évènement dangereux. Face aux enjeux importants pour la santé, l'environnement et les biens, la sécurité fonctionnelle des SIS doit être évaluée, par exemple en accord avec la CEI 61508 (IEC, 2002). Cette norme présente une approche générique des activités liées au cycle de vie de sécurité du SIS. Des exigences qualitatives ont pour objectif d'éviter les défaillances dites « systématiques ». Pour les défaillances « aléatoires », des critères quantitatifs permettent de déterminer le niveau d'intégrité de sécurité (SIL). Par exemple, pour une fonction de sécurité faiblement sollicitée, la probabilité moyenne de défaillance dangereuse à la demande (PFDavg) i.e. l'indisponibilité moyenne de la fonction de sécurité, doit être évaluée.

En accord avec la CEI 61508 (IEC, 2002) et les hypothèses précisées dans la suite, certains paramètres doivent être pris en compte dans l'évaluation de la PFDavg : architecture du système, taux des défaillances dangereuses et non détectées, intervalles des tests de révision complets, intervalles et efficacités des tests de révision partiels. Les tests complets font référence aux essais périodiques permettant de détecter toutes les défaillances d'un SIS, de telle sorte qu'il puisse être rétabli dans des conditions « comme neuves ». Les tests partiels ne permettent quant à eux que de détecter certaines défaillances, laissant les autres non-détectées jusqu'au test complet. Des inspections visuelles, des contrôles incomplets et des essais imparfaits sont des exemples de tests partiels. Bien que moins efficaces, ils peuvent être préférés aux tests complets pour plusieurs raisons :

– les tests complets sont généralement « physiques » (e.g. stimulation des éléments sensibles d'un capteur), ce qui est couteux et requière un certain temps, ils peuvent alors parfois être substitués à des tests purement « électroniques », mais qui ne couvrent pas toutes les défaillances ;

– les tests complets impliquent souvent des arrêts de production (e.g. coupure de l'alimentation électrique, d'un flux), parfois inacceptables en termes de coûts, des tests partiels (e.g. fermeture au quart de tours d'une vanne) sont alors préférés ;

– certains systèmes de sécurité ne peuvent pas être pleinement testés sans dégradation ou destruction (e.g. disques de ruptures, murs coupe feu) ;

– un test ne peut prétendre être complet qu'en conditions réelles, hors, dans de nombreux cas cela pourrait entrainer plus de danger que de prévention (e.g. détection de flamme, de gaz toxique, de surpression).

Certaines méthodes sont mentionnées à titre informatif dans la CEI 61508 (IEC, 2002) pour évaluer la PFDavg : blocs diagrammes de fiabilité, arbres de défaillance, chaînes de Markov etc. Ces outils ont notamment été comparés et discutés par Rouvroye *et al.*, 1999 & 2002, et Bukowski, 2005.

Une analyse Markovienne est utilisée par Zhang *et al.*, 2003, et des blocs diagrammes de fiabilité par Guo *et al.*, 2007, mais sans prendre en compte les tests partiels. Les problématiques de tests partiels peuvent, dans une certaine mesure, être interprétées comme de la maintenance imparfaite. Un état des lieux des méthodes et des politiques optimales de ce type de maintenance est présenté par Pham *et al.*, 1996. Par exemple, certains modèles utilisent une probabilité constante pour qu'une action de maintenance soit parfaite ou imparfaite (Brown *et al.*, 1983, Nakagawa *et al.*, 1987). Plus spécifiquement, la probabilité pour qu'une défaillance ne soit pas détectée à la suite d'une inspection est prise en compte par Badía *et al.*, 2002, mais uniquement pour un système sans redondance.

Rouvroye *et al.*, 1999 & 2002, soutiennent que des analyses Markoviennes étendues permettent de prendre en compte la majeur partie des critère de sécurité d'un système. Les modèles Markoviens sont aussi encouragés par Bukowski, 2005. Cependant, parce que les tests partiels et complets ont généralement lieux à des instants déterministes (e.g. de façon périodique), le processus n'est pas Markovien (Bukowski, 2005). Des modèles Markoviens étendus ont alors été développées afin de prendre en compte ce type de tests (Bukowski, 2001, Levitin *et al.*, 2006, Kumar *et al.*, 2008). Des approches similaires ont par exemple été utilisées à des fins d'optimisation des coûts (Wang *et al.*, 2003).

Dans le travail présenté ici, des expressions générales, obtenues par une approche analytique, sont proposées afin d'évaluer la PFDavg. Ces formulations permettent de répondre assez simplement et directement aux problématiques d'estimations paramétriques, d'évaluation et d'optimisation, en fonction des propriétés du système et des politiques de tests. La section suivante expose les hypothèses générales, les notations utilisées et les expressions proposées. Quelques applications des ces expressions sont ensuite présentées dans la troisième section.

## 2. Expressions des probabilités de défaillance à la demande (PFD)

### 2.1. *Hypothèses générales*

– Toutes les défaillances prises en compte dans la modélisation sont qualifiées de dangereuses et sont uniquement détectées lors de tests partiels ou complets ;

– Le système est composé de $N$ composants indépendants (i.e. les causes communes de défaillance ne sont pas prises en compte) et dont les taux de défaillance sont identiques et constants dans le temps ;

– Le système est d'architecture *MooN* (i.e. le système est capable d'accomplir sa fonction de sécurité si $M$ ou plus de ses $N$ composants sont opérants, avec $M \leq N$) ;

– Les $N$ composants du système sont opérants à l'instant $t_0$ ;

– Les tests partiels ne permettent de détecter que certaines défaillances définies ;

– Les tests complets permettent de détecter toutes les défaillances ;

– Les *N* composants du système sont testés simultanément lors de chaque test ;

– Les défaillances détectées au cours d'un test partiel ou complet sont réparées immédiatement et, pendant ce temps, des mesures sont prises afin de maintenir l'EUC dans un état sûr, de telle sorte que les durées de tests ou de maintenances peuvent être raisonnablement exclues de la quantification ;

– Après chaque test complet, le système est « aussi bon que neuf », la probabilité moyenne de défaillance à la demande (PFDavg) peut ainsi être calculée sur l'intervalle de test complet.

**2.2. *Notations***

| | |
|---|---|
| *MooN* | architecture du système, avec $M \leq N$ |
| $\lambda$ | taux de défaillance de chacun des *N* composants du système |
| $A_e(t)$ | disponibilité de chacun des *N* composants du système à l'instant *t* i.e. probabilité que le composant soit opérant à l'instant *t* |
| $A(t)$ | disponibilité du système à l'instant *t* i.e. probabilité que le système soit capable d'accomplir sa fonction de sécurité à l'instant *t* |
| $U(t)$ | indisponibilité du système à l'instant *t* i.e. $U(t) = 1-A(t)$ |
| $t_i$ | instant d'occurrence du $i^{ème}$ test (qui peut être partiel ou complet), avec la condition initiale : $t_0 = 0$ |
| $T_i$ | intervalle de temps entre le $(i-1)^{ème}$ et le $i^{ème}$ test i.e. $T_i = t_i - t_{i-1}$ |
| *E* | efficacité des tests partiels i.e. les défaillances détectables par chacun des tests partiels correspondent à une proportion *E* du taux de défaillance de chacun des composants |
| *n* | nombre total de tests dans un intervalle de test complet i.e. *(n-1)* tests partiels, plus le $n^{ème}$ test qui est un test complet |
| $\tau$ | intervalle de test complet i.e. $\tau = t_n$ |
| $PFD_i$ | probabilité moyenne de défaillance du système à accomplir sa fonction de sécurité à la demande i.e. indisponibilité moyenne du système, entre le $(i-1)^{ème}$ et le $i^{ème}$ test i.e. sur l'intervalle $[t_{i-1}, t_i]$ |
| *PFDavg* | probabilité moyenne de défaillance du système à accomplir sa fonction de sécurité à la demande i.e. indisponibilité moyenne du système, sur l'intervalle de test complet i.e. $[0, \tau]$ |
| $Obs_i$ | probabilité d'observer la défaillance de chacun des composants du système lors du $i^{ème}$ test |
| $k_i$ | nombre de composants défaillants observés lors du $i^{ème}$ test |
| *K* | nombre total de composants observés lors de chacun des tests |

Les notations sont reportées sur la Figure 1. Les propriétés d'un système peuvent être décrites par l'ensemble *{M, N, λ}*, et la politique de tests de façon équivalente par l'ensemble *{E, $t_1$, $t_2$, ..., $t_n$}* ou *{E, $T_1$, $T_2$, ..., $T_n$}*.

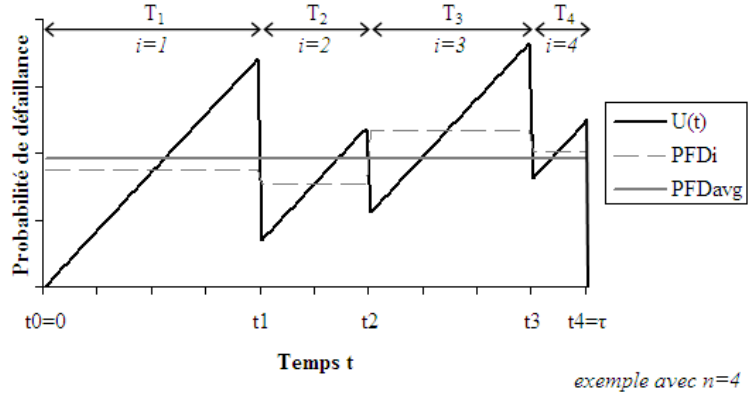

**Figure 1.** *Notations*

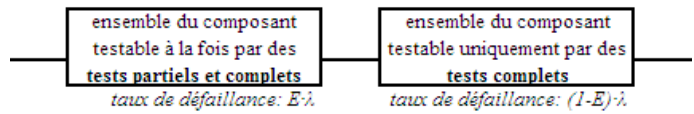

**Figure 2.** *Bloc diagramme de fiabilité pour chacun des composants du système*

### 2.3. *Expressions générales*

Le bloc diagramme de fiabilité de chacun des composants du système est représenté sur la Figure 2. En accord avec ce dernier, la disponibilité de chacun des *N* composants du système à l'instant *t* est donné par :

$$A_e(t) = e^{-E \cdot \lambda \cdot (t - t_{i-1})} \cdot e^{-(1-E)\lambda \cdot t} = e^{E \cdot \lambda \cdot t_{i-1}} \cdot e^{-\lambda \cdot t} \quad \text{pour} \quad t \in [t_{i-1}, t_i[ \quad [1]$$

La disponibilité du système d'architecture *MooN* à l'instant *t* est alors la suivante (voir démonstration en annexe) :

$$A(t) = \sum_{x=M}^{N} S(M, N, x) \cdot e^{x \cdot E \cdot \lambda \cdot t_{i-1}} \cdot e^{-x \cdot \lambda \cdot t} \quad \text{pour} \quad t \in [t_{i-1}, t_i[ \quad [2]$$

Avec la somme suivante, indépendante du temps *t* :

$$S(M, N, x) = \sum_{k=M}^{x} \binom{N}{x} \cdot \binom{x}{k} \cdot (-1)^{x-k} \quad [3]$$

Des valeurs de la somme *S(M,N,x)*, pour différentes architectures *MooN*, sont données dans le Tableau 1.

|       | $N = 1$ | $N = 2$  | $N = 3$    | $N = 4$        |
|-------|---------|----------|------------|----------------|
| $M = 1$ | 1       | 2 / -1   | 3 / -3 / 1 | 4 / -6 / 4 / -1 |
| $M = 2$ |         | 1        | 3 / -2     | 6 / -8 / 3     |
| $M = 3$ |         |          | 1          | 4 / -3         |
| $M = 4$ |         |          |            | 1              |

**Tableau 1.** *Valeurs de S(M,N,M) / S(M,N,M+1) / ... / S(M,N,N)*

La probabilité moyenne de défaillance du système à accomplir sa fonction de sécurité à la demande, sur l'intervalle *[$t_{i-1}$, $t_i$]*, est alors :

$$PFD_i = 1 - \sum_{x=M}^{N} S(M,N,x) \cdot e^{-x \cdot (1-E) \lambda \cdot t_{i-1}} \cdot \frac{1 - e^{-x \cdot \lambda \cdot T_i}}{x \cdot \lambda \cdot T_i} \qquad [4]$$

La probabilité moyenne de défaillance du système à accomplir sa fonction de sécurité à la demande, sur l'intervalle *[0, τ]*, est quant à elle :

$$PFDavg = 1 - \sum_{x=M}^{N} S(M,N,x) \cdot \sum_{i=1}^{n} \left[ e^{-x \cdot (1-E) \lambda \cdot t_{i-1}} \cdot \frac{1 - e^{-x \cdot \lambda \cdot T_i}}{x \cdot \lambda \cdot \tau} \right] \qquad [5]$$

Lorsque le produit $\lambda \cdot \tau$ est suffisamment faible (i.e. $\lambda \cdot \tau \ll 10^{-2}$), les approximations suivantes peuvent être obtenues par un développement de Taylor :

$$A_e(t) \approx 1 + E \cdot \lambda \cdot t_{i-1} - \lambda \cdot t \qquad \text{pour} \quad t \in [t_{i-1}, t_i[ \qquad [6]$$

$$A(t) \approx 1 - \binom{N}{M-1} \cdot \lambda^{N-M+1} \cdot (t - E \cdot t_{i-1})^{N-M+1} \qquad \text{pour} \quad t \in [t_{i-1}, t_i[ \qquad [7]$$

$$PFD_i \approx \binom{N}{M-1} \cdot \frac{\lambda^{N-M+1}}{N-M+2} \cdot \frac{1}{T_i} \cdot \left( (t_i - E \cdot t_{i-1})^{N-M+2} - (t_{i-1} \cdot (1-E))^{N-M+2} \right) \qquad [8]$$

$$PFDavg \approx \binom{N}{M-1} \cdot \frac{\lambda^{N-M+1}}{N-M+2} \cdot \frac{1}{\tau} \cdot \sum_{i=1}^{n} \left( (t_i - E \cdot t_{i-1})^{N-M+2} - (t_{i-1} \cdot (1-E))^{N-M+2} \right) \qquad [9]$$

**2.4. *Cas particulier d'une politique de tests sans test partiel***

Sans test partiel, les expressions [2], [5], [7] et [9] deviennent :

$$A(t)^{(s)} = \sum_{x=M}^{N} S(M,N,x) \cdot e^{-x \cdot \lambda \cdot t} \approx 1 - \binom{N}{M-1} \cdot (\lambda \cdot t)^{N-M+1} \text{ pour } t \in [0, \tau[ \qquad [10]$$

$$PFDavg^{(s)} = 1 - \sum_{x=M}^{N} S(M,N,x) \cdot \frac{1 - e^{-x \cdot \lambda \cdot T_i}}{x \cdot \lambda \cdot T_i} \approx \binom{N}{M-1} \cdot \frac{(\lambda \cdot \tau)^{N-M+1}}{N-M+2} \qquad [11]$$

## 3. Applications

### 3.1. *Description du cas test*

Un système de prévention des incendies par inertage de l'atmosphère avec de l'azote sert ici d'application. Plus spécifiquement, c'est le système de mesure de la teneur en oxygène qui est étudié. Afin d'empêcher les départs et la propagation d'incendie, tout en maintenant l'atmosphère respirable, la teneur en oxygène doit être maintenu autour de 15%. L'entrepôt considéré comprend 6 capteurs d'oxygènes. Comme l'azote introduit se répartit rapidement et de façon homogène, les capteurs sont supposés redondants. La détection d'un seuil haut ou d'un seuil bas de la teneur en oxygène est alors effectuée selon une architecture en *2oo6*.

D'après les prescriptions du fabricant, chaque capteur doit être testé annuellement avec un contrôle des résultats de mesure suivi, si nécessaire, de réajustements. Ces tests sont donc supposés complets. De plus, des inspections visuelles de façon occasionnelle sont conseillées, avec quelques vérifications électroniques permises par des voyants lumineux. Il s'agit ainsi de tests partiels qui ne permettent de détecter que certaines défaillances visibles de l'extérieures. Pour des raisons de coûts fixes dus aux procédures de tests, l'ensemble des six capteurs sont testés à chaque test partiel ou complet. La politique de tests de base consiste à effectuer les tests partiels périodiquement tous les 3 mois. Au bout du $12^{ème}$ mois, le test partiel est inclus dans le test complet (i.e. *n = 4* et $T_1 = T_2 = T_3 = T_4 = 3$ *mois*, soit $t_1 = 3$ *mois*, $t_2 = 6$ *mois*, $t_3 = 9$ *mois* et $t_4 = \tau = 12$ *mois*).

L'objectif de la démarche proposée est d'utiliser les observations faites au cours des tests (i.e. les trois tests partiels effectués au $3^{ème}$, $6^{ème}$ et $9^{ème}$ mois, et le test complet effectué au $12^{ème}$ mois) afin d'estimer le taux de défaillance des capteurs et l'efficacité des tests partiels, puis d'en déduire les probabilités de défaillance à la demande (PFD). Dans un second temps, une optimisation de la répartition des tests partiels sera proposée, permettant de réduire la PFDavg.

### 3.2. *Estimations paramétriques et évaluation des PFD*

À chaque $i^{ème}$ test de chacun des capteurs, il y a une probabilité $Obs_i$ d'observer la défaillance du capteur. D'après le bloc diagramme de fiabilité de la Figure 2 :

$$Obs_i \approx E \cdot \lambda \cdot T_i \qquad \text{pour} \quad i = 1,\ldots,(n-1) \quad [12]$$

$$Obs_i \approx E \cdot \lambda \cdot T_i + (1-E) \cdot \lambda \cdot \tau \qquad \text{pour} \quad i = n \quad [13]$$

De plus, un estimateur empirique de *$Obs_i$*, noté *$\hat{Obs}_i$*, peut être ainsi défini :

$$\hat{Obs}_i = k_i / K \qquad \text{pour} \quad i = 1,\ldots,n \quad [14]$$

Les estimateurs suivants du taux de défaillance $\lambda$ des capteurs, et de l'efficacité $E$ des tests partiels, sont alors déduits des expressions [12] à [14] :

$$\hat{\lambda} = \frac{1}{K \cdot \tau} \cdot \left( \sum_{i=1}^{n} k_i \right) \qquad [15]$$

$$\hat{E} = \frac{\tau}{t_{n-1}} \cdot \left( \sum_{i=1}^{n-1} k_i \bigg/ \sum_{i=1}^{n} k_i \right) \qquad [16]$$

À noter que les observations $k_i$ suivent, par nature, une loi binomiale. Il est donc possible de définir un intervalle de confiance pour les estimations de $\lambda$ et de $E$ par l'intermédiaire d'une loi de Fisher.

Pour cette application, 4 entrepôts, chacun comprenant 6 capteurs, ont été observés sur une durée de 4 ans. Le nombre total équivalent de capteurs observés au cours de chacun des tests est donc $K = 4 \cdot 6 \cdot 4 = 96$. Les observations faites lors des tests partiels et complets ont permis d'obtenir les résultats suivants : $k_1 + k_2 + k_3 = 16$ et $k_4 = 35$. D'après [15] et [16], on en déduit alors une estimation du taux de défaillance $\lambda$ de chacun des capteurs de *6,1·10⁻⁵ heure⁻¹*, et une estimation de l'efficacité $E$ des tests partiels de *0.42*. D'après [5], la PFDavg du système d'architecture *2oo6* de mesure d'oxygène, en accord avec la politique de tests de base (i.e. tests partiels périodiques), est évaluée à *2,06·10⁻³*. Les PFD du système sont représentées sur la Figure 3.

### 3.3. *Optimisation de la répartition des tests partiels*

Une optimisation de la politique de tests va ici consister à répartir les instants d'occurrence des tests partiels de façon à réduire la PFDavg. Lorsque les coûts liés aux tests partiels sont indépendants de l'instant où ils sont effectués, cette démarche peut ainsi permettre de réduire la PFDavg sans surcoût. Les instants d'occurrence optimaux de chacun des tests partiels, notés $t_i^*$ avec $i = 1, ..., n-1$, (on note de même $T_i^* = t_i^* - t_{i-1}^*$ pour $i = 1, ..., n-1$) sont ainsi obtenus par la résolution du problème d'optimisation suivant :

$$\{t_1^*, t_2^*, ..., t_{n-1}^*\} = \arg \min_{0 \leq t_1 \leq t_2 \leq ... \leq t_{n-1}} [PFDavg] \qquad [17]$$

Avec PFDavg telle qu'exprimée par [5]. Ici, 3 tests partiels sont à répartir au cours d'une année. La résolution de [17] conduit alors aux résultats suivants : $T_1^* = $ *4,8 mois*, $T_2^* = $ *3,0 mois*, $T_3^* = $ *2,3 mois*, $T_4^* = $ *1,9 mois* (i.e. $t_1^* = $ *4,8 mois*, $t_2^* = $ *7,8 mois*, $t_3^* = $ *10,1 mois* et $t_4 = \tau = $ *12 mois*) qui impliquent une PFDavg de *1,87·10⁻³*, soit une réduction d'environ 10% par rapport à la politique de tests de base. De plus, l'indisponibilité maximale du système sur l'intervalle de test complet a été réduite de plus de 25%. Les PFD du système, en accord avec cette politique de tests optimisée, sont représentées sur la Figure 4.

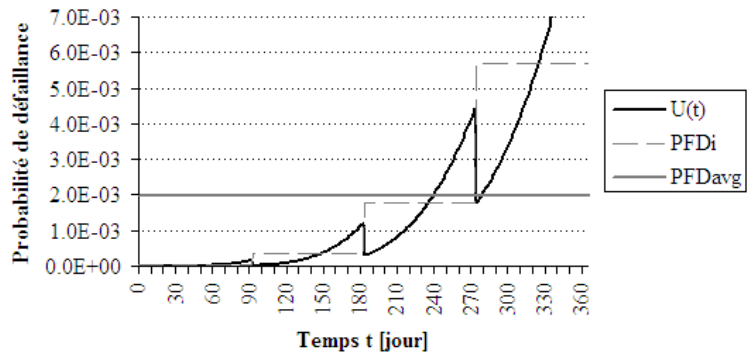

**Figure 3.** *PFD du système de mesure d'oxygène selon la politique de tests de base*

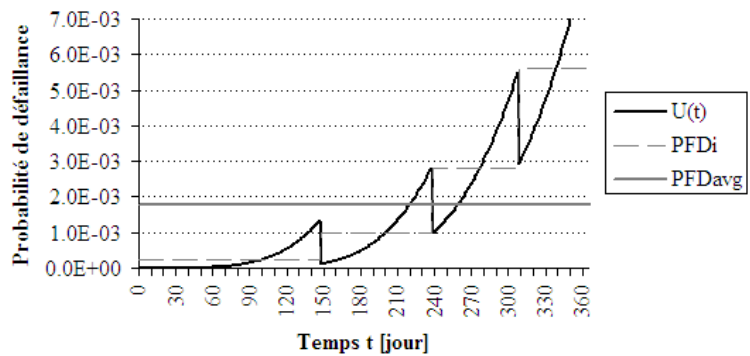

**Figure 4.** *PFD du système de mesure d'oxygène selon la politique optimisée*

## 4. Conclusion

Les expressions générales, introduites pour l'évaluation des PFD d'un système d'architecture *MooN* soumis à des tests partiels et complets, permettent de disposer d'un outil relativement simple pour le management des risques industriels. Un moyen d'estimer les performances du système et des politiques de tests a notamment été présenté, à partir du retour d'expérience obtenu lors des tests partiels et complets. De plus, il a été montré qu'une optimisation de la répartition des tests partiels permet d'améliorer certains critères de sécurité. Pour l'application proposée, une réduction d'environ 10% de la PFDavg et de plus de 25% de l'indisponibilité maximale a par exemple été permise, par rapport à une politique de tests partiels périodiques. Sous les hypothèses présentées, il est ainsi possible d'optimiser les performances du système sans surcoût.

Parmi les hypothèses utilisées, nombreuses sont celles qui peuvent être relâchées par quelques modifications des expressions présentées, par exemple en sommant des termes qui modélisent les causes communes de défaillance ou l'indisponibilité du système en phases de test ou de maintenance. D'autres perspectives concernent notamment le recours à des tests asynchrones (Rouvroye 2006) i.e. les composants du système ne sont pas testés aux mêmes instants, ce qui permet également de réduire la PFDavg à moindres coûts.

**5. Annexe : démonstrations**

Les résultats suivants sont proposés pour $t_{i-1} \leq t < t_i$ avec $i = 1, ..., n$.

D'après [1] et les hypothèses présentées au 2.1, la disponibilité du système d'architecture *MooN* à l'instant $t$ est :

$$A(t) = \sum_{k=M}^{N} \left[ \binom{N}{k} \cdot A_e(t)^k \cdot (1 - A_e(t))^{N-k} \right] = \sum_{k=M}^{N} \left[ \binom{N}{k} \cdot e^{k \cdot E \cdot \lambda \cdot t_{i-1}} \cdot e^{-k \cdot \lambda \cdot t} \cdot \left(1 - e^{E \cdot \lambda \cdot t_{i-1}} \cdot e^{-\lambda \cdot t}\right)^{N-k} \right]$$

En utilisant le binôme de Newton :

$$A(t) = \sum_{k=M}^{N} \left[ \binom{N}{k} \cdot e^{k \cdot E \cdot \lambda \cdot t_{i-1}} \cdot e^{-k \cdot \lambda \cdot t} \cdot \sum_{l=0}^{N-k} \left[ \binom{N-k}{l} \cdot (-1)^{N-k-l} \cdot \left( e^{(N-k-l) E \cdot \lambda \cdot t_{i-1}} \cdot e^{-(N-k-l)\lambda \cdot t} \right) \right] \right]$$

$$A(t) = \sum_{k=M}^{N} \sum_{l=0}^{N-k} \left[ \binom{N}{k} \cdot \binom{N-k}{l} \cdot (-1)^{N-k-l} \cdot e^{(N-l) E \cdot \lambda \cdot t_{i-1}} \cdot e^{-(N-l)\lambda \cdot t} \right]$$

D'après le théorème de Fubini, il est possible de permuter ainsi les sommes :

$$A(t) = \sum_{l=0}^{N-M} \sum_{k=M}^{N-l} \left[ \binom{N}{k} \cdot \binom{N-k}{l} \cdot (-1)^{N-k-l} \cdot e^{(N-l) E \cdot \lambda \cdot t_{i-1}} \cdot e^{-(N-l)\lambda \cdot t} \right]$$

En utilisant le changement de variable $x = N - l$ :

$$A(t) = \sum_{x=M}^{N} \sum_{k=M}^{x} \left[ \binom{N}{x} \cdot \binom{x}{k} \cdot (-1)^{x-k} \cdot e^{x \cdot E \cdot \lambda \cdot t_{i-1}} \cdot e^{-x \cdot \lambda \cdot t} \right]$$

Ce qui permet de poser :

$$A(t) = \sum_{x=M}^{N} \left[ S(M,N,x) \cdot e^{x \cdot E \cdot \lambda \cdot t_{i-1}} \cdot e^{-x \cdot \lambda \cdot t} \right] \quad \text{avec} \quad S(M,N,x) = \sum_{k=M}^{x} \left[ \binom{N}{x} \cdot \binom{x}{k} \cdot (-1)^{x-k} \right]$$

Les probabilités moyennes de défaillance à la demande sont alors :

$$PFD_i = \frac{1}{T_i} \cdot \int_{t_{i-1}}^{t_i} U(t) \cdot dt = 1 - \frac{1}{T_i} \cdot \int_{t_{i-1}}^{t_i} A(t) \cdot dt = 1 - \sum_{x=M}^{N} \left[ S(M,N,x) \cdot e^{-x \cdot (1-E)\lambda \cdot t_{i-1}} \cdot \frac{1 - e^{-x \cdot \lambda \cdot T_i}}{x \cdot \lambda \cdot T_i} \right]$$

$$PFDavg = \frac{1}{\tau} \cdot \sum_{i=1}^{n} \left[ T_i \cdot PFD_i \right] = 1 - \sum_{x=M}^{N} \left[ S(M,N,x) \cdot \sum_{i=1}^{n} \left[ e^{-x \cdot (1-E)\lambda \cdot t_{i-1}} \cdot \frac{1 - e^{-x \cdot \lambda \cdot T_i}}{x \cdot \lambda \cdot \tau} \right] \right]$$

**6. Bibliographie**